\documentstyle[12pt]{article}
\newcommand{\const}{\mathop{\rm const}\limits}

\newcommand{\vraisup}{\mathop{\rm vraisup}\limits}

\newcommand{\supp}{\mathop{\rm supp}\limits}

\begin{document}

\begin{center}

{\bf TCHEBYCHEV'S CHARACTERISTIC OF REARRANGEMENT INVARIANT SPACE.  } \\

\vspace{4mm}

   E.Ostrovsky and L.Sirota, {\sc Israel}. \\

\vspace{2mm}

{\it Department of Mathematics and Statistics, Bar-Ilan University,
59200, Ramat Gan.}\\
e \ - \ mails: eugostrovsky@list.ru;  \ galo@list.ru; \ sirota@zahav.net.il \\

\vspace{3mm}
                     {\sc Abstract.}\\

\end{center}

\vspace{3mm}

  We introduce and investigate in this short article a new characteristic of rearrangement invariant
(r.i.) (symmetric) space, namely  so-called  Tchebychev's characteristic. \par
  We reveal an important class of the r.i. spaces - so called regular r. i. spaces and show that the
 majority of known r.i. spaces: Lebesgue-Riesz, Grand Lebesgue Spaces, Orlicz, Lorentz and Marcinkiewicz r.i.
 spaces are regular. But we construct after several examples of r.i. spaces without the regular property. \par

 Applications - Probability theory and Statistics.\par

\vspace{3mm}

  \par

\vspace{2mm}

\par

\vspace{3mm}

\section{Notations. Statement of problem.}
\vspace{3mm}

 Let $ ( \Omega, \cal{A}, \mu )  $ be measure space with sigma-finite non-trivial measure $ \mu $
and $ (F, ||\cdot|| =   ||\cdot||F )  $  be any rearrangement invariant (r.i.) space over
$ ( \Omega, \cal{A}, \mu ).  $ \par

 The detail investigate of r.i. spaces see in the classical books \cite{Bennet1}, \cite{Krein1}. \par

  Hereafter $ C, C_j $ will denote any non-essential finite positive constants. As usually, for the measurable function
 $ f: \Omega \to R, $

  $$
|f|_p(\Omega, \mu) = |f|_p(\mu) = |f|_p = \left[\int_{\Omega} |f(\omega)|^p \ \mu(d \omega) \right]^{1/p}, \ 1 \le p < \infty,
 $$
  $ L_p(\Omega,\mu) = L_p(\mu) = \{f: \ |f|_p < \infty \}; \ m $
 will denote usually Lebesgue measure, and we will write in this case $ m(dx) = dx; $
  $ |f|_{\infty} \stackrel{def}{=} \vraisup_{\omega} |f(\omega)|. $ \par

 We will conclude without loss of generality in the case when $  \mu(\Omega) < \infty \
\Rightarrow \mu(\Omega) = 1, $   call: "probabilistic case"  and denote $ {\bf P} = \mu,  $ \par

$$
{\bf E} \xi = \int_{\Omega} \xi(\omega) \ {\bf P}(d \omega).
$$

 We presume for example construction that the source measurable space $ ( \Omega, \cal{A}, \mu )  $  is sufficiently
rich; it is suffices to set $ \Omega = [0,1]  $   or $ \Omega = [0,\infty) $ with Lebesgue measure $ m. $ \par

 We denote as usually for arbitrary finite a.e.measurable function (random variable) $ \xi(\omega) $ its
{\it Tail function } $ T_{\xi}(t) $ as follows:

$$
 T_{\xi}(t) = \mu \{ \omega: |\xi(\omega)  \ge t  \}, \ t > 0.
$$
 The left inverse function to the  $ T_{\xi}(t)  $ is denoted $ \xi^*(t). $
\vspace{3mm}

{\bf Definition 1.1.}  Tchebychev's characteristic $ T_F(t), t > 0  $ of the {\it space} $ (F, ||\cdot|| =   ||\cdot||F )  $
is defined as follows:

$$
T^{(F)}(t) = T^{(F, ||\cdot||)}(t) \stackrel{def}{=} \sup_{\xi: \xi \in F, ||\xi||F = 1} T_{\xi}(t). \eqno(1.1)
$$
\vspace{3mm}

 {\bf Our aim is to investigate the function  $ T_F(t)  $  for sufficiently greatest values} $ t: t > t_0 = \const > 0 $
{\bf for different classes of r.i. spaces} $ (F, ||\cdot|| =  ||\cdot||F ).  $  \par

\vspace{3mm}

A possible applications of tail estimates: Functional Analysis, see the  classical books of
C.Bennet and R.Sharpley \cite{Bennet1},  S.G. Krein Yu.V. Petunin and E.M. Semenov \cite{Krein1},
 and  also  in the articles \cite{Ostrovsky22},  \cite{Ostrovsky33}, \cite{Ostrovsky44}; Probability Theory
\cite{Ostrovsky55}, \cite{Ostrovsky77}; Numerical Methods Monte-Carlo \cite{Ostrovsky66}; Statistics \cite{Ostrovsky1},
\cite{Ostrovsky2}, \cite{Ostrovsky3}, theory of random processes and fields \cite{Kozachenko1}, \cite{Ostrovsky4} etc. \par

 For  instance, let $ \theta_n  $ be $ w(n), \ w(n) \to \infty  $  at $ n \to \infty, \ n  $ is volume of sample, consistent
statistical estimate of an unknown parameter $ \theta $ for  which

$$
|| \ w(n) |\theta_n - \theta| \ ||F \le \sigma.
$$
  We can construct the confidence interval for the value $  \theta $ by means of inequality

 $$
 {\bf P}( w(n) |\theta_n - \theta| \ge u) \le T^{(F)}(u/\sigma).
 $$

\vspace{3mm}

\section{Simple properties of  Tchebychev's characteristic. Examples.}

\vspace{3mm}

 {\bf A.} Note that

$$
T^{(F)}(t) = \sup_{\xi: \xi \in F, ||\xi||F \le 1} T_{\xi}(t). \eqno(2.1)
$$
 Moreover,

$$
 \sup_{\xi: \xi \in F, ||\xi||F = C} T_{\xi}(t) =
 \sup_{\xi: \xi \in F, ||\xi||F \le C} T_{\xi}(t) = T^{(F)}(t/C), \ C = \const > 0. \eqno(2.2)
$$

\vspace{3mm}

{\bf B.} Let on the space $ (F,||\cdot||)  $ be an other norm $ |||\cdot|||   $ for which

$$
||\xi|| \ge C_1 \ |||\xi|||,  \ 0 < C_1 = \const < \infty.
$$
 Then
  $$
T^{(F, ||\cdot||)}(t)  \le T^{(F, |||\cdot|||)}(t/C_1). \eqno(2.3)
  $$
 Analogously, if

$$
||\xi|| \le C_2 \ |||\xi|||,  \ 0 < C_2 = \const < \infty,
$$
then

  $$
T^{(F, ||\cdot||)}(t)  \ge T^{(F, |||\cdot|||)}(t/C_2). \eqno(2.4)
  $$

\vspace{3mm}

{\bf C.} {\bf Definition 2.1.}   A two tail functions $ T_1(t) $ and $ T_2(t) $  are equivalent, write:
$ T_1(\cdot) \sim T_2(\cdot) $ iff there exist three
finite positive constants $ t_0 > 0, \ 0 < C_1 \le C_2 < \infty $ for which

$$
T_2(t/C_2) \le  T_1(t) \le T_2(t/C_1), \ t \ge t_0. \eqno(2.5)
$$

 We will write also $  T_1(\cdot) << T_2(\cdot),   $ iff

$$
  T_1(t) \le T_2(t/C_1), \ t \ge t_0. \eqno(2.6)
$$
 Evidently, the relation $  "<<" $ is partial order on the set of all tail functions an the
relation $  "\sim" $  is the  relation of equivalence. Also if $ T_1 << T_2 $ and $ T_2 << T_1, $
then $ T_1 \sim T_2. $ \par

\vspace{3mm}

{\bf Corollary 2.1.} If two norms on the space $ F \ ||\cdot|| $  and $ |||\cdot||| $ are equivalent in the usually
sense, then

  $$
T^{(F, ||\cdot||)}(\cdot)  \sim T^{(F, |||\cdot|||)}(\cdot). \eqno(2.7)
  $$

 As we will see further, the converse proposition is'nt true.\par
\vspace{3mm}

 Recall that the measure space  is said to be resonant, if it is non-atomic  or conversely
 completely atomic with all the atoms having equal  measure, see \cite{Bennet1}, chapter 2, section 7.\\

{\bf D. Theorem 2.1.} Let the measure $ \mu(\cdot)  $ be resonant; then for any r.i. space $  (F, ||\cdot||) $

$$
T^{(F, ||\cdot||)}(t) \le \frac{C_3(F)}{t}. \eqno(2.8)
$$

{\bf Proof.} It is known, see \cite{Bennet1}, chapter 2, section 2, theorem 2.2 that in the considered case

$$
||\xi||F \ge C_4 |\xi|_1.
$$
 We use further Tchebychev's inequality:

 $$
 T^{(L_1)}(t) \le C_5/t, \ t > 0.
 $$
 The assertion of the theorem 2.1 follows from the inequality (2.3). \par

\vspace{3mm}

{\bf E. Examples.} 1. Classical Lebesgue-Riesz spaces.\\
 We have in the case $ \mu(\Omega) = \infty $ and $ 1 \le p < \infty:  $

$$
T^{(L_p)}(t) = t^{-p}, \ t > t_0.
$$
 When $ \mu(\Omega) = 1 $

$$
T^{(L_p)}(t) = \min \left(1, t^{-p} \right), \ t > t_0. \eqno(2.9)
$$
 Indeed, the upper estimate it follows from Tchebychev's inequality; the lower estimate  follows
from the consideration of following example:

$$
{\bf P} (\xi=t) = t^{-p}; \  {\bf P} (\xi=0) = 1-t^{-p}, \ t > 1.
$$

 Obviously,

 $$
 T^{(L_{\infty})}(t) = 0, \ t > 1.
 $$

\vspace{3mm}

2. Generalized Lorentz space.\\
 Let $ \mu = {\bf P} $ and let
 $ w = w(t), \ t > 0  $ be positive continuous strictly  increasing function, $ \lim_{t \to \infty} w(t) = \infty.  $
A  generalized Lorentz space $ L^{(w)} $ consists by definition on all the measurable functions $ \xi(\omega) $
with finite norm (more precisely, quasinorm)

$$
||\xi||L^{(w)}  = \sup_{t > 0} \left[ w(t) \ T_{\xi}(t) \right]. \eqno(2.10)
$$

 We conclude as before

 $$
 T^{L^{(w)}}(t) = \min (1,1/ w(t) ). \eqno(2.11)
 $$

 {\bf Remark 2.1.} We observe if $ w(t) = t^{p}, $ then $ T^{L^{(w)}}(t) = T^{L_p}(t), \ t > 1,  $
but the spaces  $ L^{(w)} $ and $ L_p $ are not isomorphic.\par

\vspace{3mm}

\section{Tchebychev's characteristic and fundamental functions. Regular r.i. spaces.}
\vspace{3mm}

{\it We study in this section the relations between Tchebychev's characteristic and fundamental functions.} \par

 We impose on the measure $  \mu $ here the restriction that it is diffuse: for arbitrary measurable set $ B $
there is its measurable subset $  D $ such that

$$
\mu(D) =  \mu(B)/2.
$$

   Recall that a fundamental function  $  \phi_F(\delta), \ \delta \in (0,\infty) $
  of the r.i. space $  (F, ||\cdot||)  $ over the measurable space
$ ( \Omega, \cal{A}, \mu )  $ may be defined as follows:

$$
\phi_F(\delta)
\stackrel{def}{=} \sup_{D:  \mu(D) \le \delta} ||I(D)||F. \eqno(3.1)
$$
 Here and further $  I(D)= I(D,\omega) $ is an indicator function of the measurable set $ D. $ \par
The application of  fundamental function in the functional analysis, in particular, in the theory of
interpolation of operators is described  in \cite{Bennet1}, \cite{Krein1}; the application in the theory
of  approximation see in  \cite{Ostrovsky7}. \par

 Many examples of fundamental functions for different r.i. spaces are computed in the books
\cite{Bennet1}, \cite{Krein1}.  For the so-called Grand Lebesgue Spaces the fundamental functions
are investigated and calculated  in \cite{Liflyand1}, \cite{Ostrovsky22}. \par

 Let us consider for instance the case of Orlicz's space $ Or(N) $ over our measurable space, in
which we assume the measure $  \mu $ to be diffuse. We suppose also that
 the Young function $ N=N(u) $ is in addition strictly monotonic on the positive semi-axes  and continuous.\par
 We will use in this article the Luxemburg  norm in the space $ Or(N): $

$$
||\xi||Or(N) = \inf \left\{k > 0, \int_{\Omega} N \left(\frac{|\xi(\omega)|}{k} \right) \ \mu(d \omega) \le 1 \right\}.
\eqno(3.2)
$$

 The fundamental function of this space has a view

 $$
 \phi_{Or(N)}(\delta) = \frac{1}{N^{-1}(1/\delta)}. \eqno(3.3)
 $$
  Hereafter $ g^{-1}(t) $ will denote the inverse function to the function $ g(\cdot). $ \par
  Further, let $ \xi \ge 0, \ ||\xi||Or(N) = 1. $ Since  the Young function $ N=N(u) $
is strictly monotonic and continuous

$$
\int_{\Omega} N(\xi(\omega)) \ \mu(d \omega) = 1,
$$
therefore

$$
T_{\xi}(t) \le 1/N(t).
$$
 We conclude analogously to the Lebesgue-Riesz and Lorentz spaces considering the example

 $$
 {\bf P} (\xi_0 = t) = 1/N(t) = 1 - {\bf P}(\xi_0 = 0), \  t = \const: \ N(t) > 1,
 $$
for which

$$
{\bf E}N(\xi_0) = 1,
$$
 that

$$
T^{(Or(N))}(t) = 1/N(t), \ t > t_0. \eqno(3.4)
$$

\vspace{4mm}

{\bf Definition 3.1.} The r.i. space $ (F, ||\cdot||F)  $ is said to be {\it regular } r.i. space, if

$$
\left[ \frac{1}{\phi_F(1/t)} \right]^{-1} = \frac{1}{T^{(F)}(t)}. \ t > t_0. \eqno(3.5)
$$

The r.i. space $ (F, ||\cdot||F)  $ is said to be {\it weak regular } r.i. space, if

$$
\left[ \frac{1}{\phi_F(1/t)} \right]^{-1} \asymp \frac{1}{T^{(F)}(t)}. \ t > t_0. \eqno(3.6)
$$

 We have proved the following fact.\\

\vspace{3mm}

{\bf Theorem 3.1.} The Orlicz's space $ Or(N) $ over our measurable space, in
which we assume the measure $  \mu $ to be diffuse and suppose also that
the Young function $ N=N(u) $ is in addition strictly monotonic increase and continuous,
is  regular r.i. space. \par
 If we replace the Luxemburg norm on some equivalent, we obtain the weak regular space.\par
\vspace{3mm}

 {\bf Examples 3.1.} For the spaces $ L_p $ over diffuse sigma-finite measure we have

 $$
 \phi_{L_p}(\delta) = \delta^{1/p}, \ T^{(L_p)} (1/\delta) = \delta^p = \left[\phi_{L_p}(\delta) \right]^{-1}.
 $$

 Another examples of weak regular r.i. spaces are the classical Lorentz and Marcinkiewicz spaces.\par

\vspace{3mm}

\section{Tchebychev's characteristic of associate regular r.i. spaces.}
\vspace{3mm}

 Recall that the associate r.i. space $ (F', ||\cdot|| F' )    $ to the space $ (F, ||\cdot|| F ) $
consists on all the measurable functions $ g: \Omega \to R $ with finite norm

$$
||g||F' = \sup_{\xi: ||\xi||F = 1} \left| \int_{\Omega} g(\omega) \ \xi(\omega) \ \mu(d \omega) \right|. \eqno(4.0)
$$

 Under some additional conditions (absolutely continuous norm etc.) the associate space may  coincides with
conjugate  (dual) space $ (F^*, ||\cdot|| F^* ); $ for instance,  it is true for Orlicz's space
$ (\Omega, N(u)) $ iff the Young function $ N(u) $ satisfies the $ \Delta_2 $ condition. \par

\vspace{3mm}

{\bf Theorem 4.1.} Assume again that $ (\Omega, \cal{A}, \mu)  $ is resonant measure space.
Suppose also both the r.i. spaces $ (F, ||\cdot|| F ), \  (F', ||\cdot|| F' ) $
are regular.  Then

$$
\left[ \frac{1}{T^{(F)}} \right]^{-1}(t) \cdot  \left[ \frac{1}{T^{(F')} } \right]^{-1}(t) = t, \ t > 0. \eqno(4.1)
$$

\vspace{3mm}
{\bf Proof.}  Since both the r.i. spaces $ (F, ||\cdot|| F ), \  (F', ||\cdot|| F' ) $ are regular,

$$
\phi_F(\delta) = \left[ \frac{1}{T^{(F)}} \right]^{-1} \left(\frac{1}{\delta} \right), \
\phi_{F'}(\delta) = \left[ \frac{1}{T^{(F')}} \right]^{-1} \left(\frac{1}{\delta} \right). \eqno(4.2)
$$

 We will use the known identity \cite{Bennet1}, chapter 2, section 5:

 $$
 \phi_F(\delta) \cdot \phi_{F'}(\delta) = \delta.  \eqno(4.3)
 $$

 It remains to substitute in equality (4.3) expressions (4.2) and write  $ t $ instead $ 1/\delta. $\par
\vspace{3mm}
{\bf Corollary 4.1.} If $ F' = F^*, $ then

$$
\left[ \frac{1}{T^{(F)}} \right]^{-1}(t) \cdot  \left[ \frac{1}{T^{(F^*)} } \right]^{-1}(t) = t, \ t > 0. \eqno(4.4)
$$
\vspace{3mm}
{\bf Corollary 4.2.} If  both the r.i. spaces $ (F, ||\cdot|| F ), \  (F', ||\cdot|| F' ) $
are weakly  regular, then

$$
\left[ \frac{1}{T^{(F)}} \right]^{-1}(t) \cdot  \left[ \frac{1}{T^{(F')} } \right]^{-1}(t) \asymp t, \ t > 0. \eqno(4.5)
$$

\vspace{3mm}

{\bf Corollary 4.3.} The condition of theorem 4.1 is satisfied if for example the space $ F $ is Orlicz space
with continuous strictly increasing Young function $ N = N(u), \ u \ge 0.$\par

\vspace{3mm}

{\bf Corollary 4.4.} Without the condition of resonance we can guarantee only the inequality

$$
\left[ \frac{1}{T^{(F)}} \right]^{-1}(t) \cdot  \left[ \frac{1}{T^{(F')} } \right]^{-1}(t) \ge t, \ t > 0. \eqno(4.6)
$$

 This fact follows immediately from the inequality

$$
 \phi_F(\delta) \cdot \phi_{F'}(\delta) \ge \delta,  \eqno(4.7)
 $$
 see also \cite{Bennet1}, chapter 2, section 5.\\

\vspace{3mm}

\section{Tchebychev's characteristic of the direct sum of r.i. spaces.}
\vspace{3mm}

{\bf Definition 5.1.}  We define for two tail functions $ T_1(\cdot), \ T_2(\cdot) $ the following operation:

$$
T_1 \vee T_2 (t) \stackrel{def}{=} \inf_{x \in [0,1]}  \left[ T_1(t x) + T_2(t (1-x))   \right]. \eqno(5.0)
$$
 Evidently, $ T_1 \vee T_2 (t) $ is again the tail function and $ T_1  \vee T_2 (t) =  T_2 \vee T_1 (t). $ \par
 Let the r.i. spaces $ (F, ||\cdot||F ) $  and $ (G, ||\cdot||G ) $  over our measurable space have Tchebychev's
characteristic functions correspondingly $ T^{(F)}(t), \ T^{(G)}(t). $  Let also a third space $  H $ be a
(direct) sum of this spaces: $ H = F + G. $\\

 \vspace{3mm}
{\bf Theorem 5.1.}

 $$
\max \left( T^{(F)}(t), T^{(G)}(t) \right) \le  T^{(H)}(t) \le T^{(F)}(t) \vee T^{(G)}(t). \eqno(5.1)
 $$
 \vspace{3mm}
 {\bf Proof.}  The left-hand side of bilateral inequality (5.1) is proved very simple. Let $ f_0 $ be
a function (depending on the variable $ t) $  from the space $  F $  such that $ ||f||_0 F = 1, \ T_{f_0}(t) = T^{(F)}(t).  $
 Then we have for the function $ h_0 = f_0 + 0 \in H: \  ||h_0||H  = 1 $ and $ T_{h_0}(t) = T^{(F)}(t),   $ therefore

$$
 T^{(H)}(t) \ge  T^{(F)}(t)
$$
 and analogously

$$
 T^{(H)}(t) \ge  T^{(G)}(t).
$$

  We will prove now the right-hand inequality in (5.1). \par
 Let $ h: \Omega \to R $ be any function from the space $ H $  with unit norm in this space.
We can suppose without loss of generality by virtue of definition of sum of two spaces  that  exist  two
functions say  $ f,  f \in F $ and $ g, g \in G $  for which $ h = f + g $   and
$$
1 = ||h||H = ||f||F + ||g||G. \eqno(5.2)
$$
 It follows from the equality (5.2) that

 $$
 ||f||F \le 1, \  ||g||G \le 1
 $$
 and therefore

 $$
 T_{f}(t) \le T^{(F)}(t), \   T_{g}(t) \le T^{(G)}(t). \eqno(5.3)
 $$
 Let $ x $  be arbitrary number from the set $ [0,1] $ and $ y = 1-x. $ We have:

 $$
T_h(t) \le T_{f}(t x) + T_{g}(t y) \le T^{(F)}(t x) + T^{(G)}(t y).
 $$
Since the value $ x $ is arbitrary in the closed interval $ [0,1], $ we conclude

$$
T_h(t) \le  \inf_{x \in [0,1]} [ T^{(F)}(t x) + T^{(G)}(t (1-x)) ] =
[T^{(F)} \vee  T^{(G)}](t). \eqno(5.4)
$$
 Taking the supremum over $ h: ||h||H = 1 $ we obtain

 $$
 T^{(H)}(t) \le  [T^{(F)} \vee  T^{(G)}](t). \eqno(5.5)
 $$
  This completes the proof of theorem 5.1.\\

\vspace{3mm}

{\bf Example 5.1.} Let  $  F,G $ be Orlicz's  spaces over probabilistic space with diffuse measure and
with the Young functions correspondingly

$$
N_F(u) = |u|^{p_1} \ \log^{q_1}(e+|u|), \  N_G(u) = |u|^{p_2} \ \log^{q_2}(e+|u|),
$$
$ p_1,p_2 = \const > 1, q_1,q_2 = \const. $ Then the space $  H = F + G $ is also the Orlicz's space
relative the Young function $ N_h(u) = \max(N_F(u), N_G(u)  $ and
with the Tchebychev's function

$$
  T^{(H)}(t) \asymp \max \left( T^{(F)}(t), T^{(G)}(t) \right), \ t > 1.
$$

 \vspace{3mm}

\section{Tchebychev's characteristic of Grand Lebesgue-Riesz spaces (GLS).}
\vspace{3mm}

 We recall first of all in this section  for reader conventions some definitions and facts from the theory
of GLS spaces.\par

\vspace{3mm}

 Recently, see \cite{Fiorenza1}, \cite{Fiorenza2}, \cite{Fiorenza3}, \cite{Iwaniec1}, \cite{Iwaniec2},
 \cite{Kozachenko1},\cite{Liflyand1}, \cite{Ostrovsky1}, \cite{Ostrovsky2},   etc.
 appears the so-called Grand Lebesgue Spaces $ GLS = G(\psi) =G\psi =
 G(\psi; A,B), \ A,B = \const, A \ge 1, A < B \le \infty, $ spaces consisting
 on all the measurable functions $ f: X \to R $ with finite norms

     $$
     ||f||G(\psi) \stackrel{def}{=} \sup_{p \in (A,B)} \left[ |f|_p /\psi(p) \right]. \eqno(6.1)
     $$

      Here $ \psi(\cdot) $ is some continuous positive on the {\it open} interval
    $ (A,B) $ function such that

     $$
     \inf_{p \in (A,B)} \psi(p) > 0, \ \psi(p) = \infty, \ p \notin (A,B). \eqno(6.2)
     $$
We will denote
$$
 \supp (\psi) \stackrel{def}{=} (A,B) = \{p: \psi(p) < \infty, \}
$$

The set of all $ \psi $  functions with support $ \supp (\psi)= (A,B) $ will be
denoted by $ \Psi(A,B). $ \par
  This spaces are rearrangement invariant, see \cite{Bennet1}, and
  are used, for example, in the theory of probability  \cite{Kozachenko1},
  \cite{Ostrovsky1}, \cite{Ostrovsky2}; theory of Partial Differential Equations \cite{Fiorenza2},
  \cite{Iwaniec2};  functional analysis \cite{Fiorenza3}, \cite{Iwaniec1},  \cite{Liflyand1},
  \cite{Ostrovsky2}; theory of Fourier series \cite{Ostrovsky1},
  theory of martingales \cite{Ostrovsky2},mathematical statistics  \cite{Sirota2},
 \cite{Sirota4};   theory of approximation \cite{Ostrovsky7}   etc.\par

 Notice that in the case when $ \psi(\cdot) \in \Psi(A,\infty)  $ and a function
 $ p \to p \cdot \log \psi(p) $ is convex,  then the space
$ G\psi $ coincides with some {\it exponential} Orlicz space. \par
 Conversely, if $ B < \infty, $ then the space $ G\psi(A,B) $ does  not coincides with
 the classical rearrangement invariant spaces: Orlicz, Lorentz, Marcinkiewicz  etc.\par

\vspace{3mm}

{\bf Remark 6.1} If we introduce the {\it discontinuous} function

$$
\psi_r(p) = 1, \ p = r; \psi_r(p) = \infty, \ p \ne r, \ p,r \in (A,B)
$$
and define formally  $ C/\infty = 0, \ C = \const \in R^1, $ then  the norm
in the space $ G(\psi_r) $ coincides with the $ L_r $ norm:

$$
||f||G(\psi_r) = |f|_r.
$$
Thus, the Grand Lebesgue Spaces are direct generalization of the
classical exponential Orlicz's spaces and Lebesgue spaces $ L_r. $ \par

\vspace{3mm}

{\bf Remark 6.2}  The function $ \psi(\cdot) $ may be generated as follows. Let $ \xi = \xi(x)$
be some measurable function: $ \xi: X \to R $ such that $ \exists  (A,B):
1 \le A < B \le \infty, \ \forall p \in (A,B) \ |\xi|_p < \infty. $ Then we can
choose

$$
\psi(p) = \psi_{\xi}(p) = |\xi|_p.
$$

 Analogously let $ \xi(t,\cdot) = \xi(t,x), t \in T, \ T $ is arbitrary set,
be some {\it family } $ F = \{ \xi(t, \cdot) \} $ of the measurable functions:
$ \forall t \in T  \ \xi(t,\cdot): X \to R $ such that
$$
 \exists  (A,B): 1 \le A < B \le \infty, \ \sup_{t \in T} \
|\xi(t, \cdot)|_p < \infty.
$$
Then we can choose

$$
\psi(p) = \psi_{F}(p) = \sup_{t \in T}|\xi(t,\cdot)|_p.
$$
The function $ \psi_F(p) $ may be called as a {\it natural function} for the family $ F. $
This method was used in the probability theory, more exactly, in
the theory of random fields, see \cite{Ostrovsky1}. \par

 More detail investigations of tail and fundamental  functions of
GLS see in \cite{Ostrovsky1}, \cite{Ostrovsky2}, \cite{Liflyand1}.\par

 {\it We consider in this section only the cases } $ \mu = {\bf P} $  and $ B < \infty. $\\

\vspace{3mm}

 An important\\
  {\bf Example 6.1. } Let $ B = \const  > 1, \ \beta = \const > 0 $ and let

$$
\psi_{B,\beta}(p) = (B - p)^{-\beta}, \ 1 \le p < B \eqno(6.3)
$$
and $ \psi_{B,\beta}(p) = \infty  $ otherwise. \par
 For instance: if $ \Omega = (0,1), {\bf P} = m $ and $ \xi_2(\omega) = \omega^{-1/2},  $ then
 $ \xi_2(\cdot) \in G\psi_{2,1/2}(\cdot). $\par
  Notice that for all positive values $ \epsilon < 0.5 $

 $$
 \xi_2(\cdot) \notin G\psi_{2+ \epsilon,1/2}(\cdot) \cup  G\psi_{2,1/2-\epsilon}(\cdot)
 $$
and that the function $ \psi_{2,1/2}(p) $ is equivalent to the natural function for the random variable
$ \xi_2(\cdot) . $

\vspace{3mm}

{\bf Lemma 6.1.} Denote

$$
\tilde{\psi} = p \cdot \log \psi(p), \ p \in [1,B).\eqno(6.4)
$$
 Proposition:

 $$
{\bf A.} \ T^{( G(\psi)  )}(t) \le \exp \left(- \tilde{\psi}^*(\log t) \right), \ t > 2. \eqno(6.5)
 $$
where $ h^*(\cdot) $ denotes the classical Young-Fenchel, or Legendre transform:

$$
h^*(x)  = \sup_y(xy -h(y)).
$$

{\bf B.}  For the spaces $ G\psi_{B,\beta}(\cdot) $ it true also the converse inequality up to dilation:

 $$
 T^{( G(\psi_{B,\beta} ) )}(t) \ge \exp \left(- \tilde{\psi}^*(\log t/C(B,\beta)) \right), \ t > 2. \eqno(6.6)
 $$

{\bf Proof. A.} Let  $ ||\xi||G\psi = 1; $ then $ ||\xi||_p \le \psi(p), \  {\bf E} |\xi|^p \le \psi^p(p).  $
 We obtain using the Tchebychev's inequality:

 $$
 T_{\xi}(t) \le \exp  \left( (p \log t - p \log \psi(p) )  \right).
 $$
 The assertion (6.5) it follows after an optimization over $ p. $  \par
The proposition  (6.6) is proved in the article \cite{Ostrovsky2}; see also \cite{Liflyand1}. \par

\vspace{3mm}

{\bf Example 6.2.} Denote $ \psi_m(p)= p^{1/m}, \ 1 \le p < \infty, \ m = \const > 0.  $ Proposition:

$$
\xi \in G\psi_m, \ \xi \ne 0  \ \Leftrightarrow  \exists C = \const \in (0,\infty), \
T_{\xi}(t) \le \exp \left( - C \ t^m  \right).
$$

 We will  formulate the main result of this section, which may be obtained after simple calculations  basing
 on the lemma 6.1.\\

\vspace{3mm}

{\bf Theorem 6.1.} There exists a {\it non-regular} r.i. space over the probabilistic space with diffuse measure,
namely the space $ G\psi_{B,\beta}  $ with  $ B > 1, \ \beta > 0.  $  \\

{\bf Proof.} Let us consider the space $ G\psi_{B,\beta}. $
In detail:

$$
T^{(G\psi_{B,\beta})}(t) \asymp t^{-B} \ (\log t)^{ \beta B }, \ t \to \infty,
$$

$$
\phi_{ G\psi_{B,\beta} }(\delta) \asymp \delta^{1/B} \ |\log \delta|^{\beta}, \ \delta \to 0+,
$$
so that at $ t \to \infty $

$$
\left[ \frac{1}{\phi_{ G\psi_{B,\beta} }(1/t)}   \right]^{-1} \asymp t^B \ (\log t)^{ \beta B },
$$

$$
\frac{1}{T^{(G\psi_{B,\beta}) } (t) } \asymp t^B \ (\log t)^{ - \beta B }.
$$


\vspace{3mm}

\vspace{3mm}
\begin{thebibliography}{99}

\vspace{3mm}


\bibitem{Bennet1}
  Bennet C. and Sharpley R. Interpolation of operators. Orlando, Academic Press Inc., 1988.

\bibitem{Fiorenza1}
   {\sc Capone C., Fiorenza A., Krbec M.} On the Extrapolation Blowups in the
   $ L_p $ Scale. Collectanea Mathematica, {\bf 48}, 2, (1998), 71 - 88.
 \bibitem{Fiorenza2}
 {\sc Fiorenza A.} Duality and reflexivity in grand Lebesgue spaces.
       Collectanea Mathematica (electronic version), {\bf 51}, 2, (2000), 131 - 148.
\bibitem{Fiorenza3}
 {\sc Fiorenza A., and Karadzhov G.E.} Grand and small Lebesgue spaces and
       their analogs. Consiglio Nationale Delle Ricerche, Instituto per le
      Applicazioni del Calcoto Mauro Picine”, Sezione di Napoli, Rapporto tecnico n.
      272/03, (2005).
\bibitem{Iwaniec1}
   {\sc Iwaniec T., and Sbordone C.} On the integrability of the Jacobian under
      minimal hypotheses. Arch. Rat.Mech. Anal., 119, (1992), 129 – 143.
\bibitem{Iwaniec2}
 {\sc Iwaniec T., P. Koskela P., and Onninen J.} Mapping of finite distortion:
   Monotonicity and Continuity.  Invent. Math. 144 (2001), 507 - 531.
\bibitem{Kozachenko1}
 {\sc Kozachenko Yu. V., Ostrovsky E.I.} (1985). The Banach Spaces of
      random Variables of subgaussian type. {\it Theory of Probab. and Math.
      Stat.} (in Russian). Kiev, KSU, {\bf 32}, 43 - 57.
\bibitem{Krein1}
 Krein S.G., Petunin Yu.V. and  Semenov E.M. Interpolation of Linear operators. New York, AMS, 1982.
\bibitem{Liflyand1}
{\sc Liflyand E., Ostrovsky E., Sirota L.} Structural Properties of Bilateral Grand Lebesgue Spaces.
Turk. J. Math.; {\bf 34} (2010), 207-219.
\bibitem{Ostrovsky1}
 {\sc Ostrovsky E., Sirota L. }  Universal adaptive estimations and confidence  intervals in the non-parametrical statistics.
 Electronic Publications, arXiv.mathPR/0406535 v1 25 Jun 2004.
\bibitem{Ostrovsky2}
 {\sc Ostrovsky E., Zelikov Yu.}  Adaptive Optimal Nonparametric Regression and Density Estimation based on Fourier - Legendre Expansion.
 Electronic Publication, arXiv:0706.0881v1 [math.ST] 6 Jun 2007.
\bibitem{Ostrovsky3}
{\sc Ostrovsky E.,Sirota L.}  Optimal adaptive nonparametric denoising of
multidimensional-time signal. Electronic Publication,
arXiv:0809.30211v1 [physics.data-an] 17 Sep 2008.
\bibitem{Ostrovsky4}
{\sc  Ostrovsky E.I.} Exponential Estimations for Random Fields.
Moscow - Obninsk, OINPE, (1999), in Russian.
\bibitem{Ostrovsky5}
{\sc Ostrovsky E., Rogover E. and Sirota L.} Optimal Adaptive Signal Detection
and Measurement. In: Abstracts of the International Symposium on {\sc stochastic  models in reliability engineering,
life sciences and operations management,}
Beer Sheva, Israel, (2010), p. 175.

\bibitem{Ostrovsky7}
{\sc Ostrovsky E., Sirota L.} Nikol'skii-type inequalities for rearrangement invariant
spaces. arXiv:0804.2311v1 [math.FA] 15 Apr 2008.

\bibitem{Ostrovsky22}
{\sc  Ostrovsky E. and Sirota L.} Moment Banach spaces: theory and applications.
HIAT Journal of Science and Engineering, {\bf C}, Volume 4, Issues 1 - 2,
pp. 233 - 262, (2007).

\bibitem{Ostrovsky33}
{\sc  Ostrovsky E., and Sirota L.} Boundedness of operators in bilateral Grand Lebesgue spaces,
with exact and weakly exact constant calculation.
arXiv:1103.2963 [math.FA] 15 Apr 2011.

\bibitem{Ostrovsky44}
{\sc Ostrovsky E., Sirota L., Rogover E. } Integral Operators in bilateral Grand Lebesgue Spaces.
arXiv:0912.2538 [math.FA] 13 Dez 2009.

\bibitem{Ostrovsky55}
{\sc  Ostrovsky E., and Sirota L.} Tail estimates for martingale under "LLN" norming sequense.
arXiv:1207.1908v1 [math.PR] 8 Jul 2008.

\bibitem{Ostrovsky66}
{\sc  Ostrovsky E., and Sirota L.} Monte-Carlo method for multiple parametric integrals
calculation and solving of linear integral Fredholm equations of a second kind, with confidence
regions in uniform norm.
arXiv:1101.5381 v1 [math.FA] 27 Jan  2011


\bibitem{Ostrovsky77}
{\sc  Ostrovsky E., and Sirota L.} Non-improved uniform tail estimates for mormed sums of independent
random variables with heavy tails, with applications.
arXiv:1110.4879 v1 [math.PR] 21 Oct  2011


\bibitem{Sirota2}
{\sc  Ostrovsky E., and Sirota L.}  Adaptive multidimensional-time spectral Measurements
in technical diagnosis. Communications in dependability and Managements (CDQM), Vol. 9,
No 1, (2006), pp. 45-50.

\bibitem{Sirota4}
{\sc Ostrovsky E., Sirota L.} Adaptive optimal measurements in the technical diagnostics,
reliability theory and information theory. Procedings $ 5^{th} $ international
conference on the improvement of the quality, reliability and long usage of technical
systems and technological processes, (2006), Sharm el Sheikh, Egypt, p. 65-68.

\end{thebibliography}
\end{document}